\documentclass[12pt]{article}
\usepackage{enumerate}
\usepackage{hyperref}
\usepackage{amssymb}
\usepackage{amsmath}
\usepackage{graphicx}
\usepackage{fullpage}
\usepackage{amsthm}
\usepackage{float}
\newtheorem{theorem}{Theorem}
\newtheorem{lemma}{Lemma}
\newtheorem{conj}{Conjecture}
\newtheorem{remark}{Remark}
\newtheorem{proposition}{Proposition}

\newtheorem{corollary}{Corollary}

\numberwithin{equation}{section}
\numberwithin{theorem}{section}
\numberwithin{lemma}{section}
\numberwithin{conj}{section}
\numberwithin{remark}{section}
\numberwithin{proposition}{section}
\numberwithin{corollary}{section}

\DeclareMathOperator{\hook}{hook}

\title{Multi-cores, posets, and lattice paths\thanks{The authors are grateful to Adriano Garsia for warm hospitality during the first author's visit at UCSD and his role in initiating this research.}} 

\author{Tewodros Amdeberhan\thanks{ 
(tamdeber@tulane.edu) Department of Mathematics, Tulane University,
New Orleans, LA, 70118}
\and
Emily Sergel Leven\thanks{Supported by NSF grant DGE 1144086. 
(esergel@ucsd.edu) Department of Mathematics, University of California, San Diego, CA, 92093}}

\begin{document}

\maketitle

\begin{abstract}
Hooks are prominent in representation theory (of symmetric groups) and they play a role in number theory (via cranks associated to Ramanujan's congruences). A partition of a positive integer $n$ has a Young diagram representation. To each cell in the diagram there is an associated statistic called hook length, and if a number $t$ is absent from the diagram then the partition is called a $t$-core. A partition is an $(s,t)$-core if it is both an $s$- and a $t$-core. Since the work of Anderson on $(s,t)$-cores, the topic has received growing attention. This paper expands the discussion to multiple-cores. More precisely, we explore $(s,s+1,\dots,s+k)$-core partitions much in the spirit of a recent paper by Stanley and Zanello. In fact, our results exploit connections between three combinatorial objects: multi-cores, posets and lattice paths (with a novel generalization of Dyck paths). Additional results and conjectures are scattered throughout the paper. For example, one of these statements implies a curious symmetry for twin-coprime $(s,s+2)$-core partitions. 

\vspace{12pt}
\noindent \textbf{AMS subj. class. [2010]} 05A17, 05A17, 20M99

\end{abstract}


\setcounter{section}{-1}

\section{Introduction} \label{sec:intro}

\begin{figure}
\begin{center}
\includegraphics[width=1.5in]{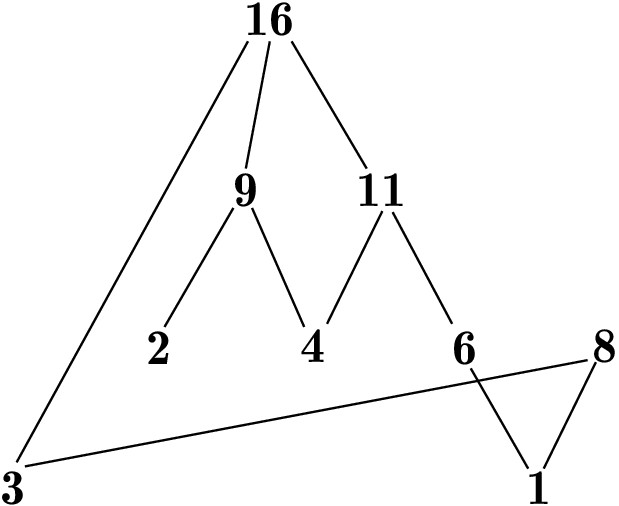}
\caption{The poset $P_{\{5,7,13\}}$}
\label{fig:Pex}
\end{center}
\end{figure}

Let \( S \) be any set of positive integers. Say that \( a \) is \emph{generated} by \( S \) if \( a \) can be written as a non-negative linear combination of the elements of \( S \).
Following the notation of \cite{Stanley}, we define \( P_S \) to be the set whose elements are positive integers not generated by \( S 
\). Equivalently, \( k \in P_S \) if \( \alpha_k = 0 \) in the generating function given by
\begin{equation}
\prod_{s \in S} \, \frac{1}{1-x^s} \, = \, \sum_{k \geq 0} \, \alpha_k \, x^k
\end{equation}
This is reminiscent of the Frobenius coin exchange problem. We make \( P_S \) into a poset by defining the cover relation so that \( a \) \emph{covers} \( b \) (written \(a \gtrdot b \)) if and only if \( a - b \in S \). For example, see Figure \ref{fig:Pex}. Note that \(P_S\) is finite if and only if the elements of \( S \) are relatively prime (no \( d > 1 \) divides every \( s \in S \)).

\begin{figure}
\begin{center}
\includegraphics[width=1.2in]{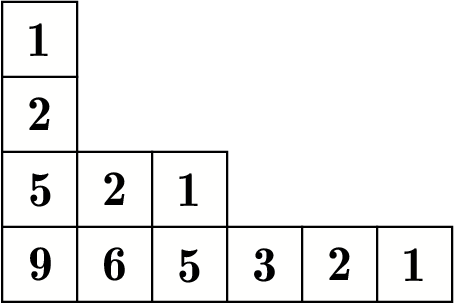}
\caption{The french Ferrers diagram of the 4-core \( (6,3,1,1) \) with hook lengths marked.}
\label{fig:hooks}
\end{center}
\end{figure}

We depict a partition \( \lambda = (\lambda_1 \geq \lambda_2 \geq \dots \geq \lambda_k > 0) \) by its French Ferrers diagram. The hook length of a cell \( c  \) in the diagram of \( \lambda \)  is the number of cells directly north or east of \( c \) including itself. It is denoted by \(\hook_\lambda(c) \), or just \(\hook(c)\) when the partition is clear. For any positive integer \(s\), we say \( \lambda \) an \(s\)-core if its diagram contains no cell \(c\) so that \( s \) divides \( \hook(c) \). For example, see Figure \ref{fig:hooks}.

Let \( (P,<_P) \) be a poset. We say that a set \( I \subseteq P\) is a \emph{lower ideal} of this poset if \( a <_P b \) and \( b \in I \) imply \( a \in I \). The work of \cite{Keith} gives a natural bijection between \(s\)-cores and lower ideals of \( P_{\{s\}} \). In particular, this bijection associates the lower ideal \( I \) of \( P_{\{s\}} \) with the \(s\)-core whose first column has hook lengths given by \( I \). For example, the 4-core in Figure \ref{fig:hooks} corresponds to the lower ideal \{1,2,5,9\} of \( P_{\{4\}} \). 

Note that this bijection will send an \(S\)-core to a set which is simultaneously a lower ideal of \( P_{\{s\}} \) for each \(s \in S\), and hence is a lower ideal of \( P_S \). For example, the lower ideal \{1,4,6,11\} of \( P_{\{5,7,13\}} \) (see Figure \ref{fig:Pex}) corresponds to the \{5,7,13\}-core partition (8,4,3,1).

When \( S = \{s, t\} \) and \( \gcd(s,t)=1\), \cite{Anderson} notes that there is a bijection between \(S\)-cores and paths in the \( s \times t \) lattice rectangle consisting of north and east steps which stay above the line \( y = \frac{t}{s} \cdot x \). We will call these paths \((s,t)\)-Dyck paths. Furthermore, \cite{Armstrong} makes multiple conjectures regarding \( \{s, t\} \)-cores. In \cite{Stanley}, the authors prove one of these conjectures in the special case \( t = s+1 \). They refer to this case as the Catalan case because the aforementioned paths are precisely Dyck paths, which are counted by the Catalan number \( C_s = \frac{1}{s+1}{2s \choose s} \).

The remainder of this paper consists of four sections. Section 1 gives a \(q\)-analog of the formula for the number of \((s,t)\)-cores. This has a nice interpretation in terms of \((s,t)\)-Dyck paths. Section 2 explores the properties of the case \( t=s+2 \). Section 3 is dedicated to developing a theory of $\{s, s+1, s+2,\dots, s+p\}$-cores. These correspond to an interesting generalization of Dyck paths, which is presented in Section 4. The authors believe this to be a new family of combinatorial objects which includes Dyck paths and Motzkin paths as special cases.


\section{Counting \((s,t)\)-Dyck paths}

As noted above, whenever \(s\) and \(t\) are coprime, there is a bijection between \((s,t)\)-cores and \((s,t)\)-Dyck paths. It is easy to show, using the Cyclic Lemma of \cite{Cyclic}, that the number of \((s,t)\)-Dyck paths is simply
\begin{equation} \label{mncount}
\frac{1}{s+t} {s+t \choose s}.
\end{equation}
However, we can get another formula using the following theorem of \cite{Krewaras}.

\begin{theorem}[\cite{Krewaras}] \label{det}
The number of partitions contained in a shape \( \lambda = (\lambda_1 \geq \lambda_2 \geq \dots \geq \lambda_k ) \) is given by
\[ \det \left( {\lambda_j + 1 \choose j-i+1} \right). \]
\end{theorem}

\begin{figure}
\begin{center}
\includegraphics[width=1in]{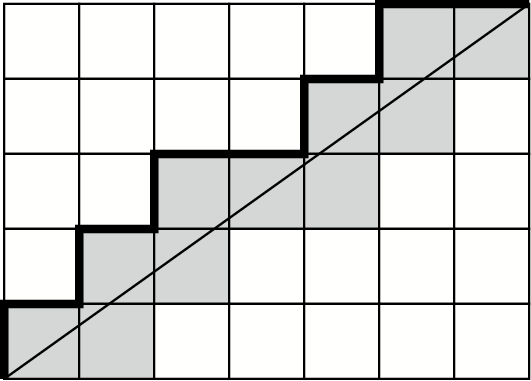}
\caption{The path closest to the main diagonal in the \(5 \times 7\) rectangle.}
\label{75hug}
\end{center}
\end{figure}

Note that the squares above each path in the \( s \times t \) rectangle form the English Ferrers diagram of a partition. Moreover, the path which hugs in the diagonal corresponds in this way to the partition with parts \( \lfloor{ \frac{s}{t} \cdot i } \rfloor \) for \(i\) from 1 to \(t-1\). For example, Figure \ref{75hug} depicts the path hugging the diagonal in the \( 7 \times 5 \) rectangle. We can see the partition \( (5,4,2,1) \) above this path.

It is surprising that in these cases the determinant of Theorem \ref{det} has the simple closed form of (\ref{mncount}). Moreover, when \(t=s+1\), (\ref{mncount}) gives the Catalan number \(C_s\) and Theorem \ref{det} yields an interesting identity.

\begin{corollary}
The Catalan numbers satisfy the identity
\[ \sum_{k=1}^{n} (-1)^k {k+1 \choose n-k} C_k = 0. \]
\end{corollary}

\begin{proof}
As observed above, we have \( C_n = \det \left( {n-j+1 \choose j-i+1} \right) \). By reversing the orders of the rows and columns (which introduces a sign of \((-1)^{2(n-1)}=1\)), we get \( C_n = \det \left( {j+1 \choose i-j+1} \right) \). The matrix here is Hessenberg. There is a recursion \cite{Muir} for the determinant of Hessenberg matrices which, in the present case, yields
\begin{align*}
C_n =  \det \left( {j+1 \choose i-j+1} \right)_{i,j=1}^{n-1} &= (-1)^{n-1} \sum_{k=1}^{n-1} (-1)^k {k+1 \choose n-k}  \det \left( {j+1 \choose i-j+1} \right)_{i,j=1}^{k-1} \cr
 &= (-1)^{n-1} \sum_{k=1}^{n-1} (-1)^k {k+1 \choose n-k} C_k. \qedhere
\end{align*}
\end{proof}

For a partition \( \lambda = (\lambda_1 \geq \lambda_2 \geq \dots \geq \lambda_k ) \), we say \(\mu\) is contained in \(\lambda\), written \( \mu \leq \lambda \), if \( \mu = ( \mu_1 \geq \mu_2 \geq \dots \geq \mu_k \geq 0) \) and \( \mu_i \leq \lambda_i \) for \(1 \leq i \leq k\). Now, we present a simple \(q\)-analog of Theorem \ref{det}.

\begin{theorem} \label{qdet}
For a partition \( \lambda = (\lambda_1 \geq \lambda_2 \geq \dots \geq \lambda_k ) \), we have the following determinantal formula:
\[
\sum_{\mu \leq \lambda} q^{|\mu|} = \det \left( q^{j-i+1 \choose 2} \left[ {\lambda_j + 1 \atop j-i+1} \right]_q \right)
\]
where $|\mu| = \mu_1 + \dots + \mu_k$ is the size of the partition $\mu$.
\end{theorem}

\begin{proof}

We will proceed by using the principle of inclusion-exclusion. Let \( A \) be the set of all sequences \(a = (a_1,a_2,\dots,a_k)\) such that \(0 \leq a_i \leq \lambda_i \) for each \(i\). We will write \( |a| = a_1 + a_2 + \dots + a_k \). Let \( A_i  \subseteq A \) be the set of such sequences with an ascent \( a_i < a_{i+1} \) at position \(i\). Then for any \( S \subseteq [k-1] \) we let \( A_S = \cap_{i \in S} \, A_i \). In particular, \( A_{\emptyset} = A\). The sequences corresponding to partitions \( \mu \leq \lambda \) are contained in \( A - \cup_{i=1}^{n-1} A_i \).

For any set \(B \subseteq A \), write
\[ q(B) = \sum_{b \in B} q^{|b|}. \]
By the principle of inclusion-exclusion,
\[
q(A - \cup_{i=1}^{n-1} A_i ) \,=\, \sum_{S \subseteq [k-1]} (-1)^{|S|} \, q(A_S).
\]
Therefore we just need to evaluate \( q(A_S) \) for each \(S\).

Note that any \( S \subseteq [k-1] \) induces a set partition of \( [k] \) as follows: \( i \) and \(i+1\) are in the same set exactly when \( i \in S \). For example,
\[ S = \{2,3,5\} \hbox{ and } n=8 \Longleftrightarrow \{ \{1\}, \{2,3,4\}, \{5,6\}, \{7\}, \{8\} \}.\]
Let this partition be denoted \( X_S \). Note that whenever \( i \) and \( j \) are in different blocks of \( X_S \), then we can choose \( a_i \) and \(a_j \) independently for \( a \in A_S \). Hence we can express \( q(A_S) \) as a product by breaking up our choices for the sequences into blocks.

Each block of \( X_S \) has the form $\{r,r+1,\dots,s\}$. We need to choose \(s-r+1\) values for the corresponding \(a_i\)'s which are strictly increasing. Since \(a_i \leq a_s \leq \lambda_s \leq \lambda_i\) for each \( r \leq i \leq s\), we may choose any values for the \(a_i\)'s of this block which satisfy \( 0 \leq a_r < a_{r+1} < \dots < a_s \leq \lambda_s \). This given, the contribution of the \(a_i\)'s in this block to \(\sum_{a} q^{|a|} \) will be 
\[
q^{s-r+1 \choose 2} \left[ {\lambda_s + 1 \atop s-r+1} \right]_q.
\]

Define \( r_i^S \) to be the smallest element of the \(i\)th block of \(X_S\) and \( s_i^S \) to be the largest element. Then we have
\[
q(A - \cup_{i=1}^{n-1} A_i ) \,=\, \sum_{S \subseteq [k-1]} (-1)^{|S|} \prod_{m=1}^{|X_S|} q^{s_m^S - r_m^S +1 \choose 2}  \left[ {\lambda_{s_m^S} + 1 \atop s_m^S - r_m^S +1} \right]_q.
\]

On the other hand, consider the matrix of Theorem \ref{qdet}. Note that the subdiagonal (where \(i=j+1\)) consists of all 1's and below that the matrix is all 0's. Therefore we can break up our expansion of the determinant
\[
\det (M_{i,j}) = \sum_{\sigma \in S_k} sign(\sigma) M_{\sigma_1,1} M_{\sigma_2,2} \dots M_{\sigma_k,k}
\]
according to which elements of the subdiagonal are chosen. Given some \( S \subseteq [k-1] \), how many \( \sigma \in S_k \) contribute to the determinant and have \( \sigma_j = j+1 \) precisely when \( j \in S \)? Exactly one. In particular, \( \sigma_{s_i^S} = r_i^S \). This is because after choosing all elements below the main diagonal, \( \sigma \) must assign the first available column \( s_1^S \) to the first available row \( r_1^S \), the next available column \( s_2^S \) to the next available row \(r_2^S\), and so on. Any other choice will include a cell below the subdiagonal, making the contribution of that permutation 0.

For example, choosing \( S=\{2,3,5\} \) when \(n=8\) induces the permutation in Figure \ref{sigchoice}. Cells from the subdiagonal which were not chosen by \(S\) and the rows and columns of chosen subdiagonals are grayed out. Then the choices within the remaining white squares are forced. We can see that in this case \( \sigma = 1 \, 3 \, 4 \, 2 \, 6 \, 5 \, 7 \, 8 \).

\begin{figure}
\begin{center}
\includegraphics[width=1.1in]{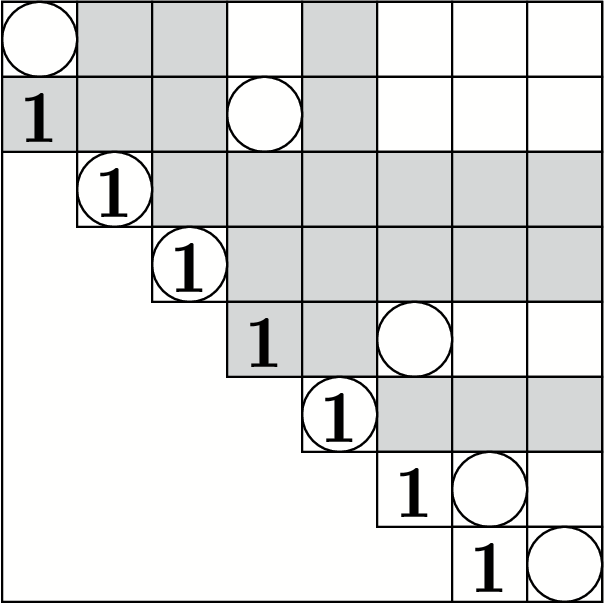}
\caption{The permutation \( 1 \, 3 \, 4 \, 2 \, 6 \, 5 \, 7 \, 8 \).}
\label{sigchoice}
\end{center}
\end{figure}

Furthermore, we have \( \sigma = (1)(2,3,4)(5,6)(7)(8) \). In general, the cycle structure of \( \sigma \) induced in this way is given by \( X_S \). Since a cycle of length \(k\) contributes \((-1)^{k-1} \) to the sign of the permutation, we have \( sign(\sigma) = (-1)^{|S|}\). Hence we obtain
\[
\det \left( q^{j-i+1 \choose 2} \left[ {\lambda_j + 1 \atop j-i+1} \right]_q \right)
\hskip -2pt = \hskip -8pt
 \sum_{S \subseteq [k-1]} \hskip -6pt (-1)^{|S|} \prod_{m=1}^{|X_S|}  q^{s_m^S - r_m^S +1 \choose 2} \left[ {\lambda_{s_m^S} + 1 \atop s_m^S - r_m^S +1} \right]_q
\hskip -6pt =
q(A - \cup_{i=1}^{n-1} A_i ). \qedhere
\]
\end{proof}

As noted above, when \(s\) and \(t\) are coprime we can use the numbers \( \lfloor \frac{s}{t} \cdot i \rfloor \) for the parts of \( \lambda \) to obtain a correspondence between \((s,t)\)-Dyck paths and partitions \( \mu \leq \lambda \). Hence this theorem gives a formula for the coarea-enumeration of \((s,t)\)-Dyck paths. In particular, this sum is
\[
\det \left( q^{j-i+1 \choose 2} \left[ { \lfloor \frac{s}{t} (t-j) \rfloor +1 \atop j-i+1} \right]_q \right).
\]

\begin{remark}
Since we have discovered this theorem, several other proofs have been found. In private correspondence, Brendon Rhoades proved the same result using the Gessel-Viennot-Lindstrom Lemma and Jeff Remmel gave a proof using involutions. Each of these proofs yielded interesting extensions.
\end{remark}


\section{A curious symmetry in \( P_{\{s,s+2\}} \) }

In order for \( P_{\{s,s+2\}} \) to be a finite poset, we must restrict ourselves to the case when \(s\) is odd. Hence we will fix \(n\) and consider \( P_{\{ 2n-1,2n+1\}} \). To sum up the results of \cite{Anderson} for this special case, we state the following theorem.

\begin{theorem}[\cite{Anderson}]
The following are equinumerous.
\begin{enumerate}
\item the number of \((s,s+2)\)-cores;
\item the number of \((s,s+2)\)-Dyck paths;
\item the number of lower ideals in \( P_{\{s,s+2\}} \).
\end{enumerate}
\end{theorem}

Consider the \( (s+1)\times(s-1) \) rectangle \(R_s\) where the (\(i,j\))-entry has value \( (s+1)(j-1) + i \), increasing North and 
Eastbound.  Abusing notation, identify $R_s$ with the set of its entries $R_s=\{1,2,\dots,(s-1)(s+1)\}$.

\begin{theorem} \label{2symm}
For each \(s \geq 3\) odd, the (\(i,j\))-entry of \(R_s\) is an element of \( P_{\{s,s+2\}} \) if and only if the (\(i,s-j\))-entry is 
not. Equivalently, for \( 1 \leq i \leq s+1 \) and \(1 \leq j \leq s-1 \),
\[
(s+1)(j-1) + i \in P_{\{s,s+2\}} \Longleftrightarrow
(s+1)(s-1-j) + i \not\in P_{\{s,s+2\}}.
\]
\end{theorem}

\begin{figure}
\begin{center}
\includegraphics[width=2.8in]{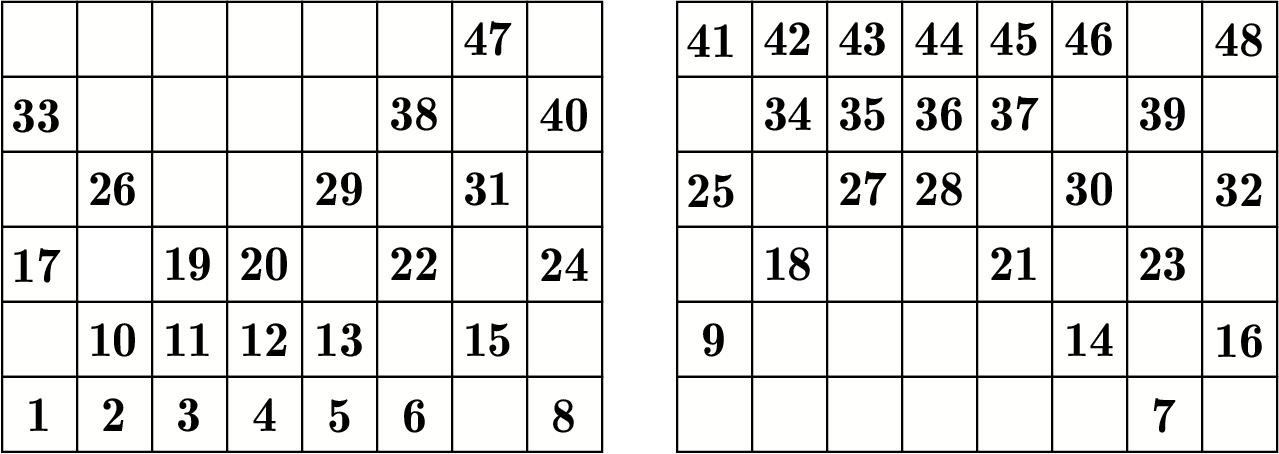}
\caption{The rectangle \(R_7\) with entries in \( P_{\{7,9\}} \) and entries not in \( P_{\{7,9\}} \).}
\label{sp2}
\end{center}
\end{figure}

Before we can prove Theorem \ref{2symm}, we need to introduce some more notation and a result of \cite{Popoviciu}. Let \( \lfloor x 
\rfloor \) and \( \{x\} \) denote the \emph{integral part} and the \emph{fractional part} of a real number \(x\) respectively. Define \( N_{s,t}(m) = 
\#\{ (k,\ell) \in \mathbb{Z}^2 \,\vert\, k,\ell \geq 0, sk+t\ell=m\} \) to be the number of partitions of \(m\) using only the elements 
\(s\) 
and \(t\) as parts. It is then clear that \[
\frac{1}{(1-x^s)(1-x^t)} \,=\, \sum_{m \geq 0} N_{s,t}(m) x^m
\]
and that $m\in P_{s,t}$ if and only if $N_{s,t}(m)=0$.

\begin{proposition}[\cite{Popoviciu}] \label{popo}
Let \( \equiv_a \) denote congruence modulo \(a\). If \(s\) and \(t\) are relatively prime, then
\[
N_{s,t}(m) \,=\, \frac{m}{st} - \left\{ \frac{t^{-1} m}{s} \right\} - \left\{ \frac{s^{-1}m}{t} \right\} + 1,
\]
where \(s^{-1},t^{-1} \in \mathbb{N}\) are relative modular inverses; that is \(t^{-1}t \equiv_s 1\) and \(s^{-1}s \equiv_t 1\).
\end{proposition}

There are two consequences of this result, the first of which is due to Sylvester.

\begin{corollary} \label{PopoCor}
If \(s\) and \(t\) are relatively prime, then
\begin{enumerate}
\item the largest integer (Frobenius number) in \(P_{\{s,t\}}\) is \(st-s-t\);
\item exactly half of the integers in \( \{1,\dots,(s-1)(t-1)\} \) belong to \(P_{\{s,t\}}\).
\end{enumerate}
\end{corollary}

\begin{proof}[Proof of Theorem \ref{2symm}]
Write \(s=2n-1\), \(t=s+2=2n+1\). Split up the set of entries in 
\(R_s= \{1,2,\dots,(s-1)(t-1)\} \) into two disjoint 
halves: 
$R_s=U \cup V$ where 
\[
U = \{ m \,\vert\, m=\alpha(s+1)+r+1, \, 0 \leq \alpha \leq \frac{s-3}{2}, \, 0 \leq r \leq s\}.
\] \[
V = \{ m' \,\vert\, m'=(s-2-\alpha)(s+1)+r+1, \, 0 \leq \alpha \leq \frac{s-3}{2}, \, 0 \leq r \leq s\},
\] 
Each $m$ and $m'$ is uniquely determined by $(\alpha,r)$. Hence the map \( \Psi: U \to V \) given by \( \Psi(m(\alpha,r))
=\Psi(\alpha(s+1)+r+1) = m'(\alpha,r)
=(s-2)(s+1)-\alpha(s+1)+r+1 \) is a bijection (where $m$ and $m'$ are given as above by the same 
values of $\alpha$ and $r$). Maintaining $m'=\Psi(m)$, introduce the subsets 
\begin{align*}
U_0&=\{m\in U \,\vert\, m\in P_{s,s+2}, m'\in P_{s,s+2}\},
   \qquad U_1=\{m\in U \,\vert\, m\in P_{s,s+2}, m'\not\in P_{s,s+2}\},\\ 
U_2&=\{m\in U \,\vert\, m\not\in P_{s,s+2}, m'\in P_{s,s+2}\},
   \qquad U_3=\{m\in U \,\vert\, m\not\in P_{s,s+2}, m'\not\in P_{s,s+2}\} 
\end{align*}
and similarly
\begin{align*}
V_0&=\{m'\in V \,\vert\, m'\in P_{s,s+2}, m\in P_{s,s+2}\},
   \qquad V_1=\{m'\in V \,\vert\, m'\in P_{s,s+2}, m\not\in P_{s,s+2}\},\\ 
V_2&=\{m'\in V \,\vert\, m'\not\in P_{s,s+2}, m\in P_{s,s+2}\},
   \qquad V_3=\{m'\in V \,\vert\, m'\not\in P_{s,s+2}, m\not\in P_{s,s+2}\}. 
\end{align*}

The key is this: \emph{exactly one} of the two integers \(m\) and \(m'\) determined by a pair \((\alpha,r)\) is in the poset \(P_{s,s+2}\). To see this, observe \(c(s+1) \in P_{\{s,s+2\}}\) if and only if \( c \in \mathbb{N} \) is odd and \(c <s\). Indeed \(c=2d\) implies \(c(s+1)=ds+d(s+2)\not\in P_{s,s+2}\). On the other hand \(c=s-2d\) implies \(c(s+1) = (s^2-2) - s - 2(d-1)(s+1)\). Note that \(s^2-2 = st-s-t \in P_{s,s+2}\) by Sylvester's result and \(P_{s,s+2}\) is closed under subtracting \(s\) or \(s+2\). Hence \((s-2d)(s+1) \in P_{s,s+2}\) whenever \(d > 0\).

Notice \(m'-m = (s-2-2\alpha)(s+1) > 0\). Therefore, we have 
that \(m'-m \in P_{\{s,s+2\}}\). That means at least one of \(m'\) and \(m\) is in \(P_{\{s,s+2\}}\). In other words, $U_3=V_3=\emptyset$. Furthermore
Corollary \ref{PopoCor} (ii) implies $\#P_{s,s+2}=\#(R_s-P_{s,s+2})$. By definition, then, $P_{s,s+2}=U_0\cup V_0\cup U_1\cup V_1$ and $R_s-P_{s,s+2}=
U_2\cup V_2$. But $\#U_1=\#V_2$ and $\#U_2=\#V_1$, so $\#(U_1\cup V_1)=\#(U_2\cup V_2)$. Hence $U_0\cup V_0$ must be empty. That means at least one of $m'$ and $m$ does \emph{not} 
belong to $P_{s,s+2}$. The key idea has been verified and the theorem is proved.
\end{proof}

\begin{remark}
A second proof of the theorem can be achieved as follows. For \( t=s+2 \) and \( m'=\Psi(m) \),  show that \( N_{s,t}(m) + N_{s,t}(m') = 1\) by using Proposition \ref{popo} and the fact that \( s^{-1} = t^{-1} = \frac{s+1}{2}=n \). The details are left to the interested reader.
\end{remark}


\section{Counting in the multi-Catalan case} \label{sec:count}

For brevity, we will write \( T_{s,p} \) for \( P_{\{s,s+1,\dots,s+p\}}\). As mentioned above, the number of lower ideals in \( T_{s,1} \) is given by the Catalan number \( C_s = \frac{1}{s+1}{2s \choose s} \). We can prove this using a simple recursion which also applies to any \( T_{s,p} \). 

\begin{lemma} \label{JrecurLemma}
Let \( J(T_{s,p}) \) be the set of all lower ideals in \( T_{s,p} \) for any \(s,p\). Then
\begin{equation} \label{Jrecur}
|J(T_{s,p})| = \sum_{i=1}^s \, |J(T_{i-p,p})| \cdot |J(T_{s-i,p})|
\end{equation}
with initial conditions \( |J(T_{\leq 0,p})|=1 \) and \( |J(T_{i,p})| = 2^{i-1} \) for \(1 \leq i \leq p\).
\end{lemma}

\begin{proof}

The initial conditions are straightforward. Following \cite{Stanley}, for \(1 \leq i \leq s-1 \), let \(J_i(T_{s,p})\) be the set of all lower ideals in \(T_{s,p}\) containing \(1,2,\dots,i-1\) but not \(i\). Also let \(J_s(T_{s,p})\) be the set of all lower ideals in \(T_{s,p}\) containing \(1,2,\dots,s-1\), so that 
\[ J(T_{s,p}) = \bigcup_{i=1}^s J_i(T_{s,p}). \]

\begin{figure}
\begin{center}
\includegraphics[width=2.2in]{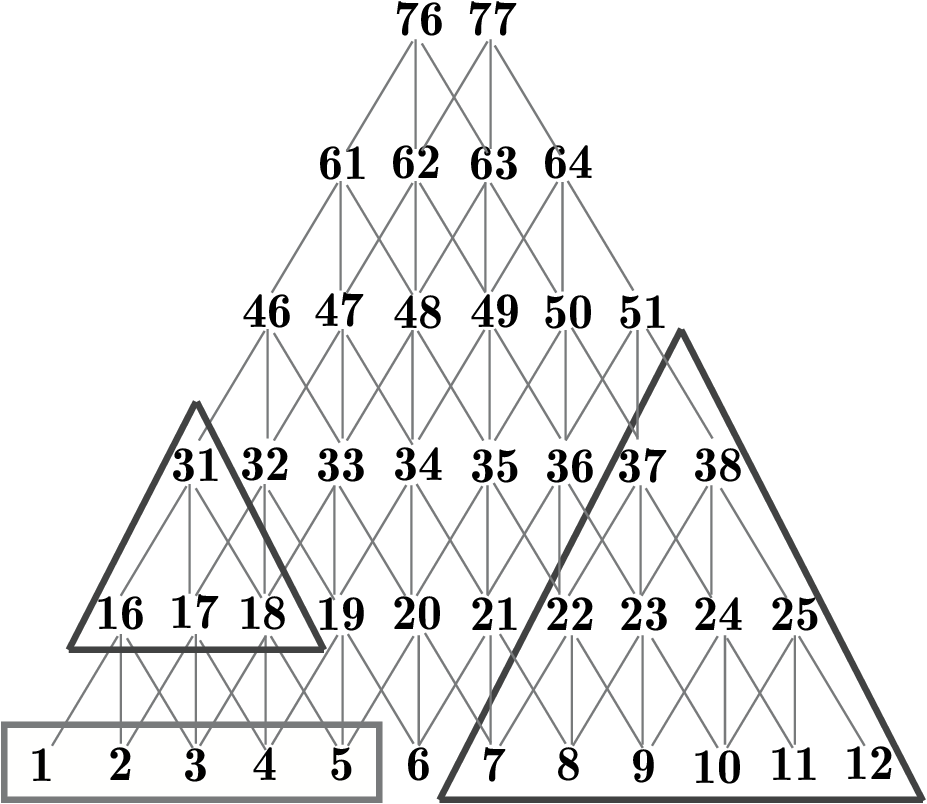}
\caption{The poset $P_{\{13,14,15\}} = T_{13,3}$. Any ideal in $J_6(T_{13,3})$ must contain the element within the rectangle and may contain elements within the triangles.}
\label{fig:J6T122}
\end{center}
\end{figure}

For example, see Figure \ref{fig:J6T122}. Here we see the poset \( T_{13,3} \) and we have marked the structure that an element of \( J_6(T_{13,2}) \) may take. Note that the structures within the triangles are isomorphic to \( T_{4,2} \) and \( T_{7,2} \), respectively. In general, the elements of \(J_i(T_{s,p})\) are in one-to-one correspondence with pairs of elements, one in \( J(T_{i-p,p}) \) and one in \( J(T_{s-i,p}) \). Hence we have the desired recursion
\begin{align*}
|J(T_{s,p})| &= \sum_{i=1}^s |J_i(T_{s,p})| \cr
&= \sum_{i=1}^s \, |J(T_{i-p,p})| \cdot |J(T_{s-i,p})|. \qedhere
\end{align*}
\end{proof}

Lemma \ref{JrecurLemma} generalizes a recursion for Catalan numbers, so we will denote \( |J(T_{s,r})| \) by \( C_s^{(r)} \). Using this recursion, we obtain the following generating function by a straight-forward computation.

\begin{theorem}
For \( r \geq 1 \) we have
\[
\sum_{s \geq 0} C_s^{(r)} x^s = \frac{ 2-2x-A_r(x) - \sqrt{ A_r(x)^2 - 4x^2 } }{2 x^{r-1}}, \hskip 12pt \hbox{ where } A_r(x) = 1 - x + \frac{x^2 - x^{r-1}}{1-x}.
\]
\end{theorem}

\noindent
As discussed above \( C_s^{(1)} \) is the classical Catalan number \( C_s \). In the next section, we will show that 
\( C_s^{(2)} \) is the \( s\)th Motzkin number. That is to say
\[
C_s^{(2)} = \sum_{k \geq 0} {s \choose 2k} C_k.
\]
Encouraged by this success, we offer the following generalized conjecture regarding partitions whose cores line up in arithmetic progression.

\begin{conj} 
Let $s$ and $d$ be two relatively prime positive integers. Then the number of \( (s,s+d,s+2d) \)-core partitions is given by
\[ \sum_{k=0}^{\lfloor\frac{s}2\rfloor} {s+d-1 \choose 2k+d-1} {2k+d \choose k}\frac1{2k+d}
=\frac1{s+d}\sum_{k=0}^{\lfloor\frac{s}2\rfloor} {s+d \choose k,k+d,s-2k}. \]
\end{conj}


\section{Generalized Dyck Paths}

Below we introduce a novel set of combinatorial objects which are in bijection with elements of \( J(T_{n,k}) \), and therefore \( (n,n+1,...n+k) \)-cores. Fix \( n, k \geq 1 \). Then \( D \) is a \emph{generalized Dyck path} if \(D\) stays above the line \(y=x\) and consists of the following:
\begin{enumerate}
\item vertical steps \(N_k = (0,k)\)
\item horizontal steps \(E_k = (k,0)\)
\item diagonal steps \(D_i = (i,i)\) for \(1 \leq i \leq k-1\).
\end{enumerate}
We will briefly denote the number of such generalized Dyck paths in the $n \times n$ square by \( GD_{n,k} \) Note that when \(k=1\), there are no diagonal steps so these paths are Dyck paths. When \(k=2\), these naturally correspond to Motzkin paths.

\begin{figure}
\begin{center}
\includegraphics[width=4.0in]{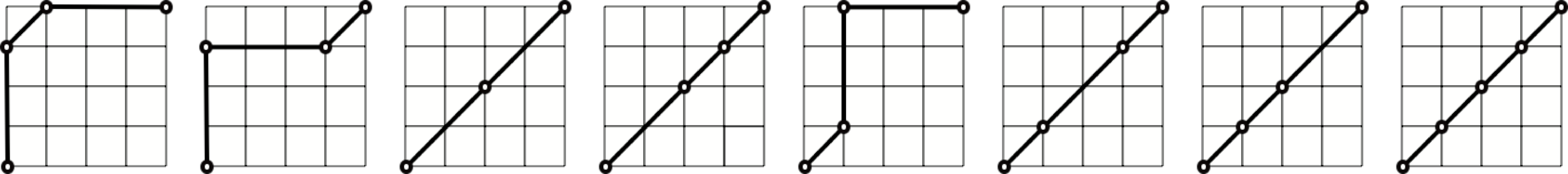}
\caption{The set of all generalized Dyck paths for \(n=4\) and \(k=3\).}
\label{GD43}
\end{center}
\end{figure}

When \( n \leq k \), we can see that \( GD_{n,k} = 2^{n-1} \) since the paths hit the main diagonal at any subset of the \( n-1 \) possible places. Moreover, by conditioning the path on the first return to the main diagonal, we get the recursion
\[ GD_{n,k} = \sum_{s=1}^{n} GD_{s-k,k} \cdot GD_{n-s,k}. \]
Here we use the convention that \( GD_{ < 0, k} = 1 \) for any \( k\). This is because when the first return is at the \(s\)th point on the main diagonal and \( s \geq k \), we must have that the first step is \( N_k \) and the last step before the return is \(E_k\). Removing these steps gives a generalized \( (s-k,k) \)-Dyck path and a generalized \( (n-s,k) \)-Dyck path. On the other hand, when \( s < k \), the first step must be \(D_s\) and the removal of this step gives a generalized \( (n-s,k) \)-Dyck path. Hence we have the following theorem.

\begin{theorem} \label{GDbij}
For any \(n,k \geq 1\), the number of generalized \((n,k)\)-Dyck paths is
\[ GD_{n,k} = C_n^{(k)}. \]
\end{theorem}

\begin{figure}
\begin{center}
\includegraphics[width=4.0in]{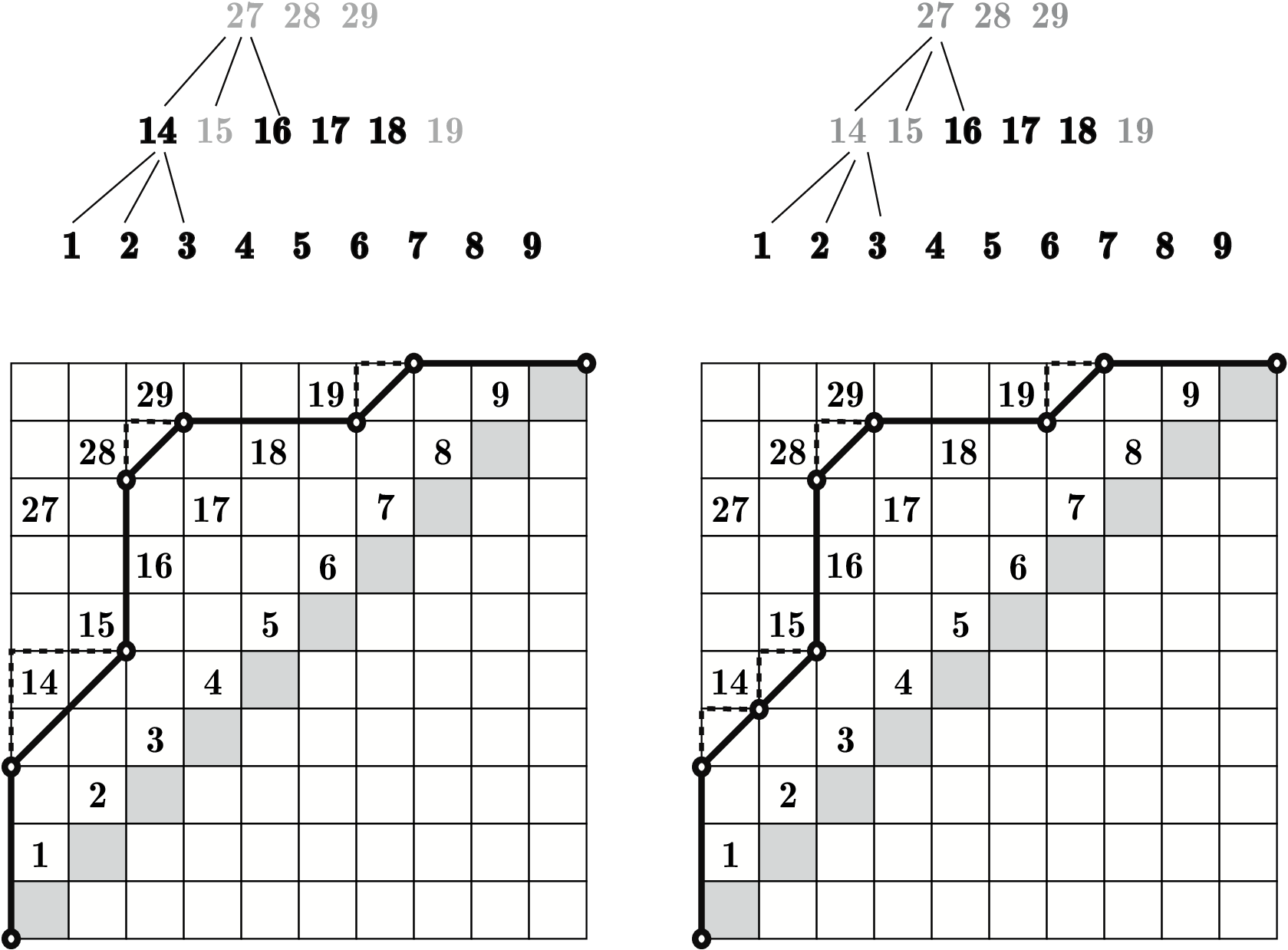}
\caption{Two lower ideals of \(T_{10,3} \) and their corresponding generalized Dyck paths. To avoid clutter, most edges of the Hasse diagrams have been excluded.}
\label{GDexample}
\end{center}
\end{figure}

\begin{proof}
In light of (\ref{Jrecur}), the recursion above is sufficient to prove the theorem. However, we can give a simple bijective proof. Simply label every \(k\)th diagonal of the \( n \times n \) lattice with consecutive numbers starting each diagonal with \(1, 1 + n+k, 1 + 2n+2k, \dots \). Then replace each step \(D_i \) by \(N_i E_i\) to get an ``inflated'' generalized Dyck path. The set of labels below the inflated path will correspond to a unique lower ideal of \(T_{n,k}\). For example, see Figure \ref{GDexample}.
\end{proof}

Since we still have the bijection between cores and lower ideals outlined in the introduction, this yields the following analog to the results of \cite{Anderson}.
\begin{theorem}
The following are equinumerous.
\begin{enumerate}
\item the number of \((s,s+1,\dots,s+k)\)-cores;
\item the number of \((s,k)\) generalized Dyck paths;
\item the number of lower ideals in \( P_{\{s,s+1,\dots,s+k\}} \).
\end{enumerate}
\end{theorem}

The largest \((a,b)\)-core was shown in \cite{Olsson} to be unique and to have size
\[
\frac{(a^2-1)(b^2-1)}{24}.
\]
Our results here enable us to construct an \((s,s+1,\dots,s+k)\)-core partition of largest size. In particular, this core corresponds to the generalized Dyck path in the \( n \times n \) lattice which takes as many vertical steps, \( N_k \), as possible followed by the longest possible diagonal step, \( D_i \) for some \(i\), and finally as many horizontal steps, \(E_k\), as possible. Looking at the image of this path under the bijection given in the proof of Theorem \ref{GDbij}, we see that the hook lengths in the first column of this core are simply
\[ 1,2,\dots,s-1, \,\, s+k+1, \dots, 2s-1, \,\, 2s+2k+1, \dots, 3s-1, \dots \]
Note that the size of a core with \( r \) parts is equal to the sum of the hooks in the first column minus \( {r \choose 2} \). This gives us an explicit but messy formula for the size of the largest \((s,s+1,\dots,s+k)\)-core. In the case \(k=2\), careful manipulation of this formula yields the following.

\begin{theorem}
The size of the largest \((s,s+1,s+2)\)-core is given by
\[
m {m+1 \choose 3} \hskip 20pt \hbox{ if }s=2m-1
\]
and
\[
(m+1) {m+1 \choose 3} + {m+2 \choose 3} \hskip 20pt \hbox{ if }s=2m.
\]
\end{theorem}

\begin{remark}
Shortly after our posting of this paper on the ArXiv, the authors of \cite{Yang} informed us (private correspondence) that they had an independent set of proofs for the results on $(s,s+1,s+2)$-core partitions.
\end{remark}

\nocite{*}
\bibliographystyle{abbrv} 
\bibliography{multicores}

\end{document}